\input amstex
\documentstyle{amsppt}
\magnification=\magstep1
\hsize=5in
\vsize=7.3in
\TagsOnRight  
\topmatter
\title On the equivalence of certain  coset conformal field theories 
\endtitle
\author Feng Xu \endauthor

\address{Department of Mathematics, 2208 Sproul Hall,
University of California at Riverside, Riverside, CA 92521}
\endaddress
\email{xufeng\@ math.ucr.edu}
\endemail
\abstract{ We demonstrate the equivalence of Kazama-Suzuki
cosets $G(m,n,k)$ and  $G(k,n,m)$ based on complex Grassmannians  by 
proving that the corresponding conformal precosheaves are isomorphic. We also
determine all the irreducible representations of
the  conformal precosheaves.
}
\endabstract

\endtopmatter
\document   
            
\heading \S1. Introduction \endheading
One of the largest two-dimensional conformal field theories (CFT) arises from
coset construction. This construction is examined from the algebraic
quantum field theory or Local Quantum Physics (LQF) (cf. [H])
point of view in [X1-3], and many mathematical 
results are obtained which have resisted  other attempts. 
Kazama and Suzuki showed in [KS] that the superconformal algebra 
based on coset $G/H$ possesses an extended $N=2$ superconformal symmetry
if, for rank$G$=rank$H$, the coset $G/H$ is a K\"{a}hler manifold. In this
paper, we forcus on the class of Kazama-Suzuki models based on the complex 
Grassmannian manifold $SU(m+n)/[SU(m)\times SU(n)\times U(1)]$. It will
be written as coset
$$           
G(m,n,k):= \frac{SU(m+n)_k \times Spin(2mn)_1}{SU(m)_{n+k} \times
SU(n)_{m+k} \times U(1)_{mn(m+n)(m+n+k)}}.
$$ 
The numerical subscripts are the levels of the representations (cf. [PS]).
The invariance of the central charge (cf. [KW]) of the coset $G(m,n,k)$
$$
c^{m,n,k}= \frac{3mnk}{m+n+k}
$$
under any permutation of $m,n,k$ suggests that the models themselves 
may be invariant [KS]. The invariance of the cosets under the exchange        
of $m$ and $n$ is manifest from their definition, but the symmetry
under the  the exchange        
of $m$ and $k$ is unexpected, as  $m$ and $k$ play rather different 
roles. In [NS], strong evidence for the symmetry is provided including
the identification of chiral quantities such as 
conformal weights, modular transformation
matrices and fusion rules under certain conditions. The goal of this
paper is to study this symmetry and related questions in the same spirit
of [X1-3]. \par 
According to the
basic idea of LQF, all the chiral quantities should be obtained by studing
the representations the conformal precosheaf of the underlying
CFT. We will recall the definition of  conformal precosheaf in \S2.
Denote by ${\Cal A}(G(m,n,k))$ the  conformal precosheaf associated
with the coset  $G(m,n,k)$. 
Hence to show that coset conformal field theory based on 
$G(m,n,k)$ is equivalent to the one based on       
$G(k,n,m)$, we just have to show that 
$${\Cal A}(G(m,n,k))\simeq {\Cal A}(G(k,n,m)) \tag 1.1
$$
where the isomorphism $\simeq$ between two
conformal precosheaves is naturally defined in \S2.1. \par
(1.1) is proved in \S3 (cf. Th. 3.7) by representing 
the two  conformal precosheaves 
on a larger Fock space and use a version of level-rank duality 
(cf. Prop. 10.6.4 of [PS]). 
An immediate corollary (cf. Cor. 3.8) is the existence of
a one to one map between the irreducible representations (primary fields)
of the two cosets and identification of all chiral quantities including
braiding and fusion matrices. However, it may be tedious to 
write down explicitly this  map in general. 
Under certain conditions, such a map is given explicitly in [NS]
which we believe to be the right one. \par
Our second goal in this paper is to determine all the irreducible
representations (primary fields )of  ${\Cal A}(G(m,n,k))$. We first
determine all the Vacuum Pairs (VPs) of the coset, a concept 
introduced in [K] which
we recall in \S2. VPs play an important role in fixed point 
resolutions and identifications of representations (cf. \S4). 
It is usually easy to come up with VPs based on
simple symmetry considerations, but it is in general a nontrival question to
determine all VPs. 
A list of VPs for $G(m,n,k)$ is given in [LVW] and [NS]  based
on Dynkin diagram symmetries, but it is known ([DJ]) that there may be
VPs which are not related to  Dynkin diagram symmetries. We show that
the list of  VPs for $G(m,n,k)$  given in [LVW] and [NS] is indeed all there
is (cf. Th. 4.4). The proof is  a mixture of solving VP equation
(2.4) for simple cases and
using the ring structure of sectors (cf. Lemma 2.7).
Using Th. 4.4 and [X3], we determine all 
the irreducible representations in Th. 4.7. \par
\par
This paper is organized as follows: In \S2.1 we give the definition of
coset conformal precosheaves and their properties. In \S.2 we recall
some basic results from [X1]  in 
Th. 2.2 and Prop. 2.3 to set up notations, and in Th. 2.4 we show that
the coset $G(m,n,k)$ has various expected properties, a result which
is implicitly contained in [X2] and [X3]. In \S2.3 we describe the notion
of Vacuum Pairs of [K] in our setting. While Lemma 2.5 follows directly
from definitions, Lemma 2.7 depends on Prop. 2.3. Lemma 2.7 plays an
important in \S4.\par
In \S3, after recalling some basic facts about the representations
of loop groups in Prop. 3.1 from [PS] and [W], we prove
Lemmas 3.1-3.6. Th. 3.7 follows from these lemmas, and Cor. 3.8 follows
from Th. 3.7. In \S4 we first recall simple selection rules about
the representations of  ${\Cal A}(G(m,n,k))$ in \S4.1. In \S4.2 we determine
all the VPs of ${\Cal A}(G(m,n,k))$ in Th. 4.4. Th. 4.4 is proved by using
Lemmas 4.1-4.3. Cor. 4.5 follows from Th. 4.4. 
Lemma 4.6 shows that the conditions of 
Lemma 2.1 of [X2] are satisfied, and so one can apply  Lemma 2.1 of [X2]
in the proof of Th. 4.7.\par  
The ideas of this paper apply to KS models based on other
Grassmannians as in [FS]. We hope to discuss those cases in the future
publication. \par
In the end of this introduction we describe in more details of the 
inclusion in the coset  $G(m,n,k)$.    
The inclusion is given by $H_1\subset G_1$
with         
$H_1=SU(m)_{n+k} \times
SU(n)_{m+k} \times U(1)_{mn(m+n)(m+n+k)}
$ and $G_1= SU(m+n)_k \times Spin(2mn)_1$. 
We will use $H_2$ and $ G_2$ to denote $H_1$ and $G_1$ respectively under the
 exchange of $m$ and $k$. 
The inclusion 
$H_1\subset G_1$ is constructed by the composition of two inclusions:
$$   
\align        
H_1  & \subset SU(m)_{n} \times SU(m)_{k}  \times SU(n)_{m} \times SU(n)_{k} \\ 
& \times U(1)_{mn(m+n)(m+n)} \times U(1)_{mn(m+n)(k)}  \tag 1.2
\endalign
$$ 
and
$$  
\align
(SU(m)_{n} & \times SU(n)_{m}  
\times U(1)_{mn(m+n)(m+n)})   \times
(SU(m)_{k} \times SU(n)_{k} \\
 &  \times U(1)_{mn(m+n)(k)})   
  \subset  Spin(2mn)_1 \times SU(m+n)_k. \tag 1.3
\endalign
$$ 
The  inclusion in (1.2) is diagonal. To describe the inclusion in (1.3), note
that the tangent space of the  Grassmanian 
$$           
\frac{SU(m+n)}{SU(m) \times
SU(n) \times U(1)}
$$  at the point corresponding to the identity of $SU(m+n)$ is isomorphic
to ${\Bbb C}^m \otimes {\Bbb C}^n$, which is a fundamental representation
of $Spin(2mn)$. The natural action of  $SU(m) \times SU(n)
\times U(1)$ on the tangent space gives the 
conformal  inclusion (cf. \S4.2 of [KW]) 
$$           
SU(m)_{n} \times SU(n)_{m} \times U(1)_{mn(m+n)(m+n)} \subset 
Spin(2mn)_1.
$$  
The inclusion 
$$
SU(m)_k \times SU(n)_k \times U(1)_{mn(m+n)(k)} \subset SU(m+n)_k
$$ comes from the conformal inclusion (cf. Prop. 4.2 of [KW])
$$
SU(m)_1 \times SU(n)_1 \times U(1)_{mn(m+n)} \subset SU(m+n)_1.
$$

\par 
\heading \S2. Preliminaries \endheading
\subheading{ \S2.1 Coset conformal precosheaf}
In this section we recall the basic properties enjoyed by the family of 
the von Neumann algebras associated with a conformal Quantum Field Theory 
on $S^1$ (cf. [GL1] ). This is an adaption of DHR analysis (cf. [H]) to
chiral CFT which is most suitable for our purposes.\par
By an {\it interval} 
 we shall always mean an open connected subset $I$
of $S^1$ such that $I$ and the interior $I^\prime $ of its complement are 
non-empty.  We shall denote by  ${\Cal I}$ the set of intervals in $S^1$.

An irreducible conformal precosheaf ${\Cal A}$ of von Neumann
algebras on the intervals of $S^1$ 
is a map
$$
I\rightarrow {\Cal A}(I)
$$
from ${\Cal I}$ to the von Neumann algebras on a Hilbert space
${\Cal H}$ that verifies the following property:
\vskip .1in
\noindent
{\bf A. Isotony}.  If $I_1$, $I_2$ are intervals and
$I_1 \subset I_2$, then
$$
{\Cal A}(I_1) \subset {\Cal A}(I_2)\, .
$$

\vskip .1in
\noindent
{\bf B. Conformal invariance}.  There is a nontrival unitary
representation $U$ of 
${\bold G}$ (the universal covering group of $PSL(2, {\bold R})$) on
${\Cal H}$ such that
$$
U(g){\Cal A}(I)U(g)^* = {\Cal A}(gI)\, , \qquad
g\in {\bold G}, \quad I\in {\Cal I} \, .
$$

The group $PSL(2, {\bold R})$ is identified with the M\"obius group of 
$S^1$, i.e. the group of conformal transformations on the complex plane 
that preserve the orientation and leave the unit circle globally 
invariant. Therefore ${\bold G}$ has a natural action on $S^1$.
\vskip .1in
\noindent
{\bf C. Positivity of the energy}.  The generator of the rotation subgroup 
$U(R)(\cdot)$ is positive.

Here $R(\vartheta )$ denotes the (lifting to ${\bold G}$ of the) rotation
by 
an angle $\vartheta $. 
\vskip .1in 
\noindent
{\bf D.  Locality}.  If $I_0$, $I$ are disjoint intervals then
${\Cal A}(I_0)$ and $A(I)$ commute.

The lattice symbol $\vee $ will denote `the von Neumann algebra generated 
by'.
\vskip .1in
\noindent
{\bf E. Existence of the vacuum}.  There exists a unit vector
$\Omega $ (vacuum vector) which is $U({\bold G})$-invariant and cyclic for
$\vee _{I\in {\Cal I}}{\Cal A}(I)$.

\vskip .1in
\noindent
{\bf F. Uniqueness of the vacuum (or irreducibility)}.  The only 
$U({\bold G})$-invariant vectors are the scalar multiples of $\Omega $.

\vskip .1in
\noindent
Assume  ${\Cal A}$ is   as defined in above.
A covariant {\it representation} $\pi $ of
${\Cal A}$ is a family of representations $\pi _I$ of the 
von Neumann algebras ${\Cal A}(I)$, $I\in {\Cal I}$, on a 
Hilbert space ${\Cal H}_\pi $ and a unitary representation
$U_\pi $ of the covering group ${\bold G}$ of $PSL(2, {\bold R})$, with 
{\it positive energy}, i.e. the generator of the rotation unitary subgroup 
has positive generator, such that the following properties hold:
$$
\align
I\supset \bar I \Rightarrow \pi _{\bar I} \mid _{{\Cal A}(I)}
= \pi _I \quad &\text{\rm (isotony)} \\
\text{\rm ad}U_\pi (g) \cdot \pi _I = \pi _{gI}\cdot 
\text{\rm ad}U(g) &\text{\rm (covariance)}\, .
\endalign
$$
A unitary equivalence class of representations of ${\Cal A}$ is called
{\it superselection sector}. \par
The composition of two superselection sectors are  known as Connes's fusion
[W]. The composition is manifestly
unitary and associative, and this is one of the most important virtues of the
above formulation. The main question is to study all 
superselection sectors of  ${\Cal A}$ and their compositions.
\par
Given two irreducible conformal precosheaves ${\Cal A}_1$, ${\Cal A}_2$
on Hilbert spaces ${\Cal H}_1$ and  ${\Cal H}_2$ with vacuum vectors
$\Omega_1$ and $\Omega_2$ respectively. One can define naturally that
 ${\Cal A}_1$ is {\it isomporphic} to  ${\Cal A}_2$ if there is a unitary map
$U: {\Cal H}_1\rightarrow  {\Cal H}_2$ such that:
$$
U^*{\Cal A}_2(I) U= {\Cal A}_1(I), \forall I\in {\Cal I};
U\Omega_1 = \Omega_2.
$$
Note that by [GL2], $U$ as defined above intertwines the representation
of the conformal group ${\bold G}$.\par 
We have the following (cf. Prop. 1.1 of [GL1]):
\proclaim{2.1 Proposition} Let ${\Cal A}$ be an irreducible 
conformal precosheaf. The
following hold:
\roster
\item"{(a)}"  Reeh-Schlieder theorem: $\Omega $ is cyclic and separating
for each von Neumann algebra ${\Cal A}(I)$, $I\in {\Cal I}$.
\item"{(b)}"  Bisognano-Wichmann property: $U$ extends to an
(anti-)unitary representation of ${\bold G}\times _{\sigma _r}{\bold Z}_2$
such that, for any $I\in {\Cal I}$,
$$
\align
U(\Lambda _I(2\pi t)) &= \Delta _I^{it} \, \\
U(r_I) &= J_I \,
\endalign
$$             
where $\Delta _I$, $J_I$ are the modular operator and the modular
conjugation associated with $({\Cal A}(I), \Omega )$.
For each $g\in {\bold G} \times _{\sigma _r} {\bold Z}_2$
$$
U(g){\Cal A}(I)U(g)^* = {\Cal A}(gI) \, .
$$
\endroster
\roster
\item"{(c)}"  Additivity: if a family of intervals $I_i$ covers the
interval $I$, then
$$
{\Cal A}(I) \subset \vee _i {\Cal A}(I_i)\, .
$$
\item"{(d)}" Haag duality: ${\Cal A}(I)' = {\Cal A}(I')$.
\endroster
\endproclaim  
Let us give some examples of conformal precosheaves.\par

Let $G$ be a compact Lie group. Denote by $LG$ the group of smooth maps
$f: S^1 \mapsto G$ under pointwise multiplication. The
diffeomorphism group of the circle $\text{\rm Diff} S^1 $ is
naturally a subgroup of $\text{\rm Aut}(LG)$ with the action given by
reparametrization. In particular ${\bold G}$ acts on $LG$. 
We will be interested
in the projective unitary representations 
(cf. Chap. 9 of [PS]) $\pi$ of $LG$ that
are both irreducible and have positive energy. This implies that $\pi $
should extend to $LG\ltimes \text{\rm Rot}$ so that
the generator of the rotation group Rot is positive.
It follows from
Chap. 9 of [PS] 
that for a fixed level there are only finite number of such
irreducible projective 
representations. \par       
Now Let $G$ be a  connected compact Lie group
and let $H\subset G$ be a Lie subgroup. 
Let $\pi^i$ be an irreducible
representations of $LG$ with positive energy at level $k$
\footnotemark\footnotetext{ When G is the direct product of simple groups,
$k$ is a multi-index, i.e., $k=(k_1,...,k_n)$, where $k_i\in \Bbb {N}$
corresponding to the level of the $i$-th simple group. The level of $LH$
is determined by the Dynkin indices of $H\subset G$. To save some writing
we write the coset simply as $H\subset G_k$ or  $H\subset G$ when the
levels are clear from the context.} 
on Hilbert space ${\bold H}^i$.
Suppose when restricting to $LH$, ${\bold H}^i$ decomposes as:
$$          
{\bold H}^i = \sum_\alpha {\bold H}_{i,\alpha} \otimes {\bold H}_\alpha, 
$$ and  $\pi_\alpha$ are irreducible representations of $LH$ on
Hilbert space ${\bold H}_\alpha$.  The set of $(i,\alpha)$ which appears in
the above decompositions will be denoted by $exp$. \par

We
shall use $\pi^1$ (resp. $\pi_1$) to denote the vacuum representation of
$LG$ (resp. $LH$) on ${\bold H^1}$ (resp. ${\bold H_1}$) . 
Let $\Omega$ (resp. $\Omega_0$)  be the vacuum vector 
in $\pi^1$ (resp.$\pi_1$) and assume 
$$
\Omega = \Omega_{0,0} \otimes \Omega_0
$$ with $\Omega_{0,0} \in {\bold H}_{1,1}$.

We shall assume that $H\subset G$ is not a conformal inclusion (cf. [KW])
to avoid triviality.\par
For each interval $I\subset S^1$, we define:
$$
{\Cal A}(I):= P \pi^1(L_IH)' \cap \pi^1(L_IG)''P,
$$  where $P$ is the projection from ${\bold H^1}$ to a closed subspace 
spanned by 
$$\vee_{J\in  {\Cal I}}\pi^1(L_JH)' \cap \pi^1(L_JG) \Omega.$$
Here $\pi^1(L_IG)''$ denotes the von Neumann algebra generated by 
$$\pi^1(a), a\in LG, \text{\rm Supp a} \subset I.$$ \par  
It follows from [X1]  that ${\Cal A}(I)$ is
an irreducible conformal precosheaf on the closed spacce.
We define this 
to be irreducible conformal precosheaf of coset ($H\subset G$) CFT and 
denote it
by ${\Cal A}_{G/H}$. 
Note the similarity of this definition to the
vertex operator algebraic definition (cf. [FZ]).  
Note also that $\pi_{(i,\alpha)}
$ above naturally gives rise to the covariant representations of 
${\Cal A}_{G/H}$. 
${\Cal A}_{G/H}$ corresponds to {\it coset} construction of
CFT. \par
For the inclusion $H_1\subset G_1$ at the end of \S1, we will also denote 
${\Cal A}_{G_1/H_1}$ by ${\Cal A}(G(m,n,k)$. \par 
\subheading{\S2.2 Some results from [X1]}
We recall some results from [X1] which will be used in the following.
We refer the reader to [X1] for more details.\par
Let $M$ be a properly infinite factor
and  $\text{\rm End}(M)$ the semigroup of
 unit preserving endomorphisms of $M$.  In this paper $M$ will always
be a type $III_1$ factor.
Let $\text{\rm Sect}(M)$ denote the quotient of $\text{\rm End}(M)$ modulo
unitary equivalence in $M$. 
It follows from
\cite{L3} and \cite{L4} that $\text{\rm Sect}(M)$  is endowed
with a natural involution $\theta \rightarrow \bar \theta $, and
$\text{\rm Sect}(M)$ is
a semiring: i.e., there are two operations $+, \times$ on 
$\text{\rm Sect}(M)$ which verifes the usual axioms. The multiplication
of  sectors is simply the composition of sectors. Hence if
$\theta_1, \theta_2$ are
two sectors, we shall write $\theta_1\times \theta_2$  as
$\theta_1\theta_2$.
In [X4],   the image of
$\theta \in \text{\rm End}(M)$ in  $\text{\rm Sect}(M)$ is denoted by 
$[\theta]$. However, since we will be mainly concerned with the ring
structure of sectors in this paper, we will denote $[\theta]$
simply by $\theta$ if no confusion arises. 
\par
Assume $\theta \in \text{\rm End}(M)$, and there exists a normal
faithful conditional expectation
$\epsilon:
M\rightarrow \theta (M)$.  We define a number $d_\epsilon$ (possibly
$\infty$) by:
$$
d_\epsilon^{-2} :=\text{\rm Max} \{ \lambda \in [0, +\infty)|
\epsilon (m_+) \geq \lambda m_+, \forall m_+ \in M_+
\}$$ (cf. [PP]).\par         
If $d_\epsilon < \infty$ for some $\epsilon$, we say $\theta$
has finite index or statistical dimension. In this case 
we define
$$
d_\theta = \text{\rm Min}_\epsilon \{ d_\epsilon |  d_\epsilon < \infty \}.
$$   $d_\theta$ is called the {\it statistical dimension} of  $\theta$. 
$d_\theta^2$ is called the {\it minimal index} of $\theta$. In fact in this
case there exists a unique $\epsilon_\theta$ such that
$ d_{\epsilon_\theta} = d_\theta$. $\epsilon_\theta$ is called the
{\it minimal conditional expectation}. 
It is clear
from the definition that  the statistical dimension  of  $\theta$ depends only
on the unitary equivalence classes  of  $\theta$.
When $N\subset M$ with $N\simeq M$, we choose $\theta \in \text{\rm End}(M)$
such that $\theta(M)=N$.   The statistical dimension 
(resp. minimal index) of the inclusion $N\subset M$ is 
defined to be the statistical dimension (resp. minimal index) of $\theta$.\par
Let  $\theta_1, \theta_2\in Sect(M)$. By Th. 5.5 of [L3],  
$d_{\theta_1+\theta_2}= d_{\theta_1} + d_{\theta_2}$, and
by Cor. 2.2 of [L5], $d_{\theta_1\theta_2}= d_{\theta_1}d_{\theta_2}$.
These two properties are usually referred to as the {\it additivity}
and {\it multiplicativity} of statistical dimensions. Also note
by Prop. 4.12 of [L4] $d_\theta = d_{\bar \theta}$. If a sector
does not have finite statistical dimension in any of the
above three equations, then the equation is understood as the 
statement that both sides of the equation are $\infty$. \par  
Assume $\lambda $, $\mu,$ and $ \nu \in \text{\rm End}(M)$ have
finite statistical dimensions.  Let
$\text{\rm Hom}(\lambda , \mu )$ denote the space of intertwiners from
$\lambda $ to $\mu $, i.e. $a\in \text{\rm Hom}(\lambda , \mu )$ iff
$a \lambda (p) = \mu (p) a $ for any $p \in M$.
$\text{\rm Hom}(\lambda , \mu )$  is a finite dimensional vector
space and we use $\langle  \lambda , \mu \rangle$ to denote
the dimension of this space.  Note that $\langle  \lambda , \mu \rangle$
depends
only on $[\lambda ]$ and $[\mu ]$. Moreover we have
$\langle \nu \lambda , \mu \rangle =
\langle \lambda , \bar \nu \mu \rangle $,
$\langle \nu \lambda , \mu \rangle    
= \langle \nu , \mu \bar \lambda \rangle $ which follows from Frobenius
duality (See \cite{L2} ).  We will also use the following
notation: if $\mu $ is a subsector of $\lambda $, we will write as
$\mu \prec \lambda $  or $\lambda \succ \mu $.  A sector
is said to be irreducible if it has only one subsector. \par 
Let $\theta_i, i=1,...,n$ be a set of irreducible 
sectors with finite index. The
ring generated by $\theta_i, i=1,...,n$ under compositions 
is defined to be a vector space 
(possibly infinite dimensional) over
${\Bbb C}$ with a  basis $\{ \xi_j , j\geq 1\}$, such that $\xi_j$ are
irreducible sectors, $\xi_j\neq \xi_{j'}$ if $j\neq j'$, and the set
$\{\xi_j, j\geq 1\}$ is a list of all   
irreducible sectors which appear as subsectors of 
finite products of $\theta_i,
i=1,...,n$. The ring 
multiplication  on the vector space is obtained naturally from that of
$Sect(M)$.\par
Let $ M(J), J\in {\Cal I}$ 
be an irreducible conformal precosheaf  on Hilbert space ${\bold H}^1$. 
Suppose
$N(J), J\in {\Cal I}$ is an irreducible conformal precosheaf 
 and $\pi^1$ is a covariant representation
of $N(J)$ on ${\bold H}^1$ such that $\pi^1(N(J)) \subset M(J)$ is a directed
standard net as defined in Definition 3.1 of [LR] for any directed
set of intervals. Fix an interval $I$ and denote by $N:= N(I), M:=M(I)$.   
For any covariant representation $\pi_\lambda$ (resp. $\pi^i$) of
the irreducible conformal precosheaf 
$N(J),   J\in {\Cal I} $ (resp. 
$M(J), J\in {\Cal I}$), let $\lambda$ (resp. $i$) be the 
corresponding  
endomorphism  of $N$ (resp. $M$) as defined in \S2.1 of [GL1]. 
These endomorphisms are obtained by localization in \S2.1 of [GL1]
and will be referred to  as  localized 
endomorphisms for convenience. The corresponding sectors will be
called localized sectors. \par
In this paper,  
if we use $1$ to denote a sector or a covariant 
representation, it should be understood as the identity sector or
vacuum representation. \par

We will use $d_\lambda$ and  $d_i$ to denote the 
statistical dimensions of $\lambda$ and  $i$ respectively. 
$d_\lambda$ and  $d_i$
are also called  the 
statistical dimensions of $\pi_\lambda$ and $\pi^i$ respectively, and they are
independent of the choice of $I$ (cf. Prop. 2.1 of [GL1]). \par
Let $\pi^i$
be a covariant representation of  $M(J),  J\in {\Cal I}$
which decomposes as:
$$
\pi^i = \sum_{\lambda} b_{i\lambda} \pi_\lambda
$$ when restricted to  $N(J),  J\in {\Cal I}$, where the sum is finite
and $ b_{i\lambda} \in {\Bbb N}$. Let $\gamma_i:=
\sum_{\lambda} b_{i\lambda} \lambda$ be the corresponding sector of $N$. 
It is shown (cf. (1) of
Prop. 2.8 in [X1]) that there are sectors $\rho, \sigma_i \in
\text{\rm Sect}(N)$ such that:
$$
\rho \sigma_i \bar \rho =\gamma_i.
$$
Notice that  $\sigma_i$ are in one-to-one correspondence with
covariant representations $\pi^i$, and in fact the map $i\rightarrow
\sigma_i$ is an isomorphism of the ring generated by $i$ and
the ring generated by $\sigma_i$.  The subfactor 
$\bar{\rho}(N) \subset N$ is conjugate to $\pi^1(N(I)) \subset M(I)$
(cf. (2) of Prop. 2.6 in [X4]). \par
Now we assume  $\pi^1(N(I)) \subset M(I)$ has finite index. Then for
each  localized sector $\lambda$ of $N$ there exists a sector
denoted by $a_\lambda$ of $N$ such that the following theorem is
true   (cf. [X4]):
\proclaim{Theorem 2.2}
(1)
The map $\lambda \rightarrow a_\lambda$ is a ring homomorphism;\par
(2)
$\rho a_\lambda = \lambda \rho, a_\lambda \bar \rho = \bar \rho \lambda
, d_\lambda = d_{a_\lambda}$;\par
(3) $\langle \rho a_\lambda,  \rho a_\mu \rangle =
\langle a_\lambda,  a_\mu \rangle = \langle  a_\lambda \bar \rho,  a_\mu
\bar \rho  \rangle$; \par
(4)  $\langle \rho a_\lambda, \rho \sigma_i  \rangle =  
 \langle a_\lambda, \sigma_i  \rangle =
\langle a_\lambda \bar \rho, \sigma_i \bar \rho
\rangle $ ; \par
(5) (3) (resp. (4)) remains valid if $ a_\lambda,  a_\mu$ (resp.  $a_\lambda$)
is replaced by any of its subsectors; \par
(6) $ a_\lambda \sigma_i =  \sigma_i a_\lambda$.
\endproclaim
We will apply the results of Th. 2.2 to  the case when 
$N(I)= {\Cal A}_{G/H}(I)\otimes
\pi_1(L_IH)''$ and $M(I)=\pi^1(L_IG)''$ under the assumption that
$H\subset G$ is {\it cofinite}, i.e., 
$\pi^1(N(I)) \subset M(I)$ has finite index, where
${\Cal A}_{G/H}(I)$ is defined in \S2.1 (cf. [X1]).  
By Th. 2.2, for every localized endomorphisms $\lambda $ of
$N(I)$ we have a map $a: \lambda \rightarrow a_\lambda$ which
verifies (1) to (6) in Th. 2.2. \par
\proclaim{Tensor Notation}
Let $\theta \in End( {\Cal A}_{G/H}(I) \otimes
\pi_1(L_IH)'')$. We will denote $\theta$ by $\rho_1\otimes \rho_2$
if 
$$
\theta(p\otimes 1)= \rho_1(p)\otimes 1, \forall p\in {\Cal A}_{G/H}(I) ,
\theta(1\otimes p')= 1\otimes \rho_2(p'), \forall p'\in\pi_1(L_IH)'' ,  
,$$ where $\rho_1 \in End( {\Cal A}_{G/H}(I) ),\rho_2 \in End(\pi_1(L_IH)'')$.
\endproclaim
Recall from \S2.1  $\pi_{i,\alpha}$ of ${\Cal A}_{G/H}(I) $ are obtained 
in the decompositions of $\pi^i$ of $LG$ with respect to subgroup $LH$,
and we  denote the set of such $(i,\alpha)$ by {\it  $exp$}. 
We will denote the sector correspondng to $\pi_{(i,\alpha)}$ simply 
by $(i,\alpha)$.  Under the condtions that $(i,\alpha), (j,\beta)$ have finite
indices,  we have that   $(i,\alpha)$ is an irreducible sector 
if and only if
$\pi_{i,\alpha}$ is an irreducible covariant representation , and
$(i,\alpha)\succ (j,\beta)$ if and only of $\pi_{j,\beta}$ appears
as a direct summand of $\pi_{i,\alpha}$,  and $(i,\alpha)$ is equal
to $(j,\beta)$ as sectors if and only  $\pi_{i,\alpha}$ is
unitarily equivalent to  $\pi_{j,\beta}$ (cf. [GL1]).\par
Given  $(i,\alpha)\in \text{\rm End}({\Cal A}_{G/H}(I)  )$ as above, we define
$(i,\alpha) \otimes 1 \in \text{\rm End} (N(I))$ so that:
$$
(i,\alpha)\otimes 1(p\otimes p') = (i,\alpha)(p)\otimes p' 
, \forall p\in {\Cal A}_{G/H}(I) , p'\in \pi_1(L_IH)''
.$$ It is
easy to see that $(i,\alpha)\otimes 1$ corresponds to the
covariant representation $\pi_{i,\alpha} \otimes \pi_1$ of
$N(I)$.
Note that this notation agrees with our tensor notation above. 
Also note that for any covariant representation $\pi_x$
of ${\Cal A}_{G/H}(I) $, we can define a localized sector $x\otimes 1$ of
$N(I)$ in the same way as in the case when $\pi_x=\pi_{i,\alpha}$. 
\par
Each covariant representation $\pi^i$ of $LG$ 
gives rise to an endomorphism
$\sigma_i\in \text{\rm End} (N(I))$ and 
$$
\rho \sigma_i \bar{\rho} = \gamma_i = \sum_{\alpha}
(i,\alpha) \otimes (\alpha) \tag 2.1
$$ where the summation is over those $\alpha$ such that $(i,\alpha) \in exp$.
The following is Prop. 4.2 of [X1]:
\proclaim{Proposition 2.3}
Assume $H\subset G$ is cofinite.  
We have:\par 
(1)
Let $x, y$ be localized 
sectors of ${\Cal A}_{G/H}(I) $ with finite index. Then
$$
\langle x, y \rangle = 
\langle a_{x\otimes 1}, a_{y \otimes 1} \rangle 
;$$\par
(2) If $(i,\alpha) \in exp$, then
$ a_{(i,\alpha)\otimes 1}  \prec a_{1\otimes \bar{\alpha}} \sigma_i$;  
\par
(3) Denote by $d_{(i,\alpha)}$ the statistical dimension of
$(i,\alpha)$. Then $d_{(i,\alpha)} \leq d_i d_\alpha$, where 
$d_i$ (resp. $d_\alpha$) is the statistical dimension of $i$ (
resp. $\alpha$). 
\endproclaim
Let us  
denote by $S_{ij}$ (resp. $\dot{S_{\alpha \beta}}$) the $S$ matrices of
$LG$        
(resp. $LH$) at level $k$ (resp. certain level of $LH$ determined by
the inclusion $H\subset G_k$) as 
defined on P. 264  of [Kac].  Define  \footnotemark\footnotetext {Our
$(j,\beta)$ corresponds to $(M,\mu)$ on P.186 of [KW], 
and $\langle (j,\beta), (1,1) \rangle$ is  equal
to $mult_M(\mu,p)$ which appears in 2.5.4 of [KW] by (2.5). So
our formula (2.2) is identical to 2.5.4 of [KW].}
$$          
b(i,\alpha) = \sum_{(j,\beta) } {S_{ij}} \overline{\dot S_{\alpha \beta}}
\langle (j,\beta), (1,1) \rangle \tag 2.2
$$          
Note the above summation is effectively over those $(j,\beta)$ such that 
$(j,\beta) \in exp$. Note that by [KW], if $b(i,\alpha)>0$, then
$(i,\alpha)\in exp$. 
The Kac-Wakimoto Conjecture (KWC) states that if
$(i,\alpha)\in exp$, then  $b(i,\alpha)>0$. Under certain conditions,
a stronger result than KWC is obtained in [X3], and the results of
[X3] apply to the coset $Gr(m,n,k)$. More precisely we have:
\proclaim{Theorem 2.4}
(1): The coset  $Gr(m,n,k)$ is cofinite 
(cf. [X1] or definition after Th. 2.2); \par
(2): There are only a finite number of irreducible representations of 
$${\Cal A}(G(m,n,k)),$$ and each  irreducible representation appears
as a direct summand of some $(i,\alpha) \in  exp$; \par
(3): The statistical dimension $d_{(i,\alpha)}$ of the coset sector
$(i,\alpha)$ 
is given by
$$
d_{(i,\alpha)} = \frac{b(i,\alpha)}{b(1,1)}
$$
where $b(i,\alpha)$ is defined in (2.2); \par
(4): The irreducible  representations of 
${\Cal A}(G(m,n,k))$ generate a unitary modular category as defined
in [Tu].
\endproclaim
\demo{Proof:}
(1) is proved at the end of \S3.2 of [X2]. (2) and (3) follows from (1)
, Cor. 3.2 and Th. 3.4 of [X3]. We note that it is assumed in 
 [X3] that all the groups involved are type $A$ groups so 
one can use the results of [W] and [X6]. But it is easy to show that
these results hold for $U(1)$ (cf. P. 58 of [X5]) since all sectors 
involved are automorphisms, and in fact it is already implicitly 
contained in \S4 of [X6]. 
Hence all results of  [X3] apply to
$U(1)$ too. (4) follows from Prop. 2.4 of [X3]. We note that
(4) also follws from (1) and [L1]. \footnotemark\footnotetext{ In 
fact using (1) and [L1] one can obtain a stronger result, i.e., 
${\Cal A}(G(m,n,k))$ is completely rational (cf. [KLM]).}
\enddemo \hfill \qed \hfill
\par  
\subheading{\S2.3 Vacuum Pairs}

Let us recall the definition of vacuum pairs according to P. 236 of [K]
(or [KW]) in our notations. As in \S2.1  
let $\pi^i$ be an irreducible
representations of $LG$ with positive energy at level $k$
on Hilbert space ${\bold H}^i$.
Suppose when restricting to $LH$, ${\bold H}^i$ decomposes as:
$$          
{\bold H}^i = \sum_\alpha {\bold H}_{i,\alpha} \otimes {\bold H}_\alpha, 
$$ and  $\pi_\alpha$ are irreducible representations of $LH$ on
Hilbert space ${\bold H}_\alpha$. 
By [GKO], the generator $L_G(0)$ of the rotation
group for  $LG$ act on ${\bold H}^i$ as
$$
L_G(0)=L_{G/H}(0)\otimes id + 1\otimes L_H(0) \tag 2.3
$$
The eigenvalues of $L_G(0)$ on ${\bold H}^i$ are given by $h_i +m, m\in 
{\Bbb Z}_{\geq 0}$, where $h_i$ is the conformal dimension or trace 
anomaly defined in 
(1.4.1) of [KW]. Let $\Omega_{i,\alpha} \otimes \Omega_{\alpha}\in 
{\bold H}_{i,\alpha} \otimes {\bold H}_\alpha$ be a unit vector 
with weight $i':=i-{\bold r}$  of 
$LG$ where  $\Omega_{\alpha}$ is the  highest weight vector of $LH$ on 
$H_\alpha$, and ${\bold r} $ is a sum of positive roots of
$LG$.   By (3.2.6) of [KW], 
$$
L_G(0) (\Omega_{i,\alpha} \otimes \Omega_{\alpha})
= (h_i+m) \Omega_{i,\alpha} \otimes \Omega_{\alpha} 
$$
where $m$ is a nonnegative integer determined by $i'$. 
But we also have
$$
\align
L_G(0) (\Omega_{i,\alpha} \otimes \Omega_{\alpha})
&= L_{G/H}(0) (\Omega_{i,\alpha})\otimes \Omega_{\alpha} + \Omega_{i,\alpha} 
\otimes L_H(0)(\Omega_{\alpha})\\
&= L_{G/H}(0) (\Omega_{i,\alpha})\otimes \Omega_{\alpha} + h_\alpha 
\Omega_{i,\alpha} \otimes \Omega_{\alpha}, 
\endalign
$$
and since the eigenvalues of $L_{G/H}(0)$ are non-negative (cf. \S3 of [KW]),
we must have
$$
h_i+m \geq  h_\alpha
$$ 
According to [K], we will say that $\{i,\alpha \}$ is a {\it Vacuum Pair} if
$$
h_i+m =  h_\alpha \tag 2.4
$$
As noted above $m$ is determined by $i':=i-{\bold r}$ and 
$\alpha$ is obtained
by restriction from weight $i'$ of $LG$ to $LH$.
Note that since there are only finitely many  $i,\alpha $, (2.4) has
only a finite number of solutions. We will denote the finite set of 
of VPs simply as  $VPS$. However it is in general nontrival to
determine $VPS$. \par 
From the equations before (2.4) we must have that (2.4) hold iff 
$$
L_{G/H}(0) (\Omega_{i,\alpha})=0,
$$
i.e.,  $\Omega_{i,\alpha}$ is a vacuum vector, and it follows immediately that
$
{\bold H}_{(1,1)}
$ 
is a direct summand of ${\bold H}_{(i,\alpha)}.$ Hence 
if the sector $(i,\alpha)$
has finite index, then  $\{i,\alpha \}$  is a VP iff
$(i,\alpha)\in exp$ and 
$$\langle (i,\alpha), (1,1) \rangle >0 \tag 2.5$$
One can see the importance of such VPs  in calculating (2.2). \par
For the rest of this section, we will assume that all sectors or
representations considered
have finite indices. \par
Assume that  
$H_1\subset H_2 \subset G.$ 
For simplicity we will use $\pi_x,\pi_y, \pi_z$ to denote the
irreducible representations of $LH_1$, $LH_2$ and $LG$ respectively,
and ${\Cal A},{\Cal B}, {\Cal C}$ to denote the conformal precosheaves of 
cosets 
$H_1\subset H_2,
H_2\subset G, H_1\subset G$ respectively. 
Note we have natural inclusions $ {\Cal A}(I)\otimes {\Cal B}(I)\subset 
{\Cal C}(I)$, 
corresponding to the
natural inclusions
$$
(\pi(L_IH_1)'\cap \pi(L_IH_2)'') \otimes
(\pi(L_IH_2)'\cap \pi(L_IG)'') \subset \pi(L_IH_1)'\cap \pi(L_IG)''
,$$
where $I$ is a proper open interval of $S^1$.
From the decompositions:
$$           
\pi_z \simeq \sum_y \pi_{(z,y)} \otimes \pi_y \simeq \sum_{y, x}
\pi_{(z,y)}  
\otimes \pi_{(y,x)} \otimes \pi_{x} 
\simeq  \sum_{x} \pi_{(z,x)} \otimes \pi_x  
$$           
we conclude that 
$$           
\pi_{(z,x)} \simeq  \sum_y \pi_{(z,y)} \otimes \pi_{(y,x)}
$$ which is understood as the decomposition of representation
$\pi_{(z,x)}$ of ${\Cal C}$ when restricted to ${\Cal A}\otimes {\Cal B} 
\subset {\Cal C}$.
The following lemma follows immediately  from (2.3) and (2.4):
\proclaim{Lemma 2.5}
$\{ z,x \}$ is a VP for $H_1\subset G$ iff there exists
a $y$ such that $\{ z,y \}$ and  $\{ y,x \}$ are VPs for $H_1\subset H_2$
and $H_2\subset G$ respectively.
\endproclaim
One can usually find VPs by examining the symmetry of Dynkin diagrams
(cf. 2.7.12  of [KW]). This motivates the following:  
\proclaim{Definition 2.6 (SVP)} 
$\{i,\alpha \} $ is called a {\it simple vacuum pair} if $d_\alpha =1$.
\endproclaim
We will denote the set of all SVPs as $SVPS$.  
\proclaim{Lemma 2.7}
SVPS is an abelian group under the compositions of sectors.
\endproclaim
\demo{Proof:} Let $\{i,\alpha \} \in SVPS$. 
By (2.5)
$\langle (i,\alpha), (1,1) \rangle >0$.
A useful property which follows from
Th. 2.2 and (2.1)   is
$$
\langle \sigma_i, a_{x\otimes \alpha}\rangle
= \langle (i,\alpha), x \rangle.
$$
Set $x=(1,1)$ we get
$\langle \sigma_i, a_{1\otimes \alpha}\rangle >0$, and so
$ \sigma_i \prec a_{1\otimes \alpha}$ since $ \sigma_i$ is
irreducible. Since $d_{ a_{1\otimes \alpha}}= d_\alpha=1$, 
$ a_{1\otimes \alpha}$ must be irreducible and
$ \sigma_i =a_{1\otimes \alpha}.$ In particular $d_i=1.$ So 
$$
a_{(i,\alpha)\otimes 1}= a_{(1,1)\otimes 1} =  
\sigma_i a_{1\otimes \bar\alpha}.
$$
It follows from  $$
a_{(\bar i,\bar \alpha)\otimes 1}=   
\sigma_{\bar i} a_{1\otimes \alpha} =a_{(1,1)\otimes 1}
$$ that  $\{\bar i,\bar \alpha \} \in SVPS$.
Now let $\{i,\alpha \}, \{j,\beta \} \in SVPS$. We must have
$ij=k, \alpha\beta = \delta$ for some
$\{\beta, \delta\}$ since all sectors have statistical 
dimension 1. To show that $SVPS$ is an abelian group, we just have to show
that $\{k,\delta \} \in VPS$. Note that
$$
a_{(k,\delta)\otimes 1}= \sigma_k a_{1\otimes \bar\delta} =
\sigma_i \sigma_j a_{1\otimes \bar\alpha} a_{1\otimes \bar\beta}
=  a_{(1,1)\otimes 1}
$$
and this shows that  $\{k,\delta \} \in VPS$ by (1) of Prop. 2.3. 
\par
\enddemo \hfill \qed \par
\heading{\S3.  ${\Cal A}(G(m,n,k)) \simeq {\Cal A}(G(k,n,m))$} \endheading
We will first recall some facts from [PS]. The reader is refered to
[PS] for more details. \par
Let ${\bold H}$ denote the 
Hilbert space $L^2(S^1; \Bbb C^N)$ of square-summable 
$\Bbb C^N$-valued functions on the circle.  The group $LU(N)$ of smooth 
maps $S^1 \rightarrow U(N)$ acts on ${\bold H}$ multiplication operators.

Let us decompose ${\bold H}= {\bold H}_+ \oplus {\bold H}_-$, where
$$
{\bold H}_+ = \{\text{\rm functions whose negative Fourier coeffients vanish}\} \, .
$$

We denote by $P $ the projection from ${\bold H}$ onto ${\bold H}_+$.

Denote by $U_{\text{\rm res}}({\bold H})$ 
the group consisting of unitary operator $A$ on ${\bold H}$
such that the commutator \cite{P, A} is a Hilbert-Schmidt operator.  
There exists a central extension $U_{\text{\rm res}}^\sim $ of
$U_{\text{\rm res}}({\bold H})$ as defined in \S 6.6 of \cite{PS}.
The central extension $ {\Cal L}U(N)$ of $LU(N)$ induced by 
 $U_{\text{\rm res}}^\sim $ is called the {\it basic extension}.

The basic representation $\pi $ of $ {\Cal L}U(N)$ is the representation 
 on Fermionic Fock space ${\Cal F}(\Bbb C^N):= \Lambda(PH)
\otimes \Lambda ((1-p )H)^*$ as defined in \S 10.6 of \cite{PS}.
Note that if $\Bbb C^N=\Bbb C^{N_1} \oplus \Bbb C^{N_2}$, then
 ${\Cal F}(\Bbb C^N)$ is canonically isomorphic to 
${\Cal F}(\Bbb C^{N_1}) \otimes  {\Cal F}(\Bbb C^{N_2}).$
\par

Let $I=\bigcup _{i=1}^n I_i$ be a proper subset of $S^1$, where $I_i$ are 
 intervals of $S^1$. Denote by $M(I, \Bbb C^N)$ the von Neumann algebra 
generated by $c(\xi )^\prime s$, with $\xi \in L^2 (I, \Bbb C^N)$.
Here $c(\xi ) = a(\xi ) + a(\xi )^*$ and $a(\xi )$ is the creation 
operator defined as in Chapter 1 of \cite{W2}.
Let $K: {\Cal F}(\Bbb C^N) \rightarrow {\Cal F}(\Bbb C^N) $ 
be the Klein transformation 
given by multiplication by 1 on even forms and by $i$ on odd forms.
We will denote the set of even forms as ${\Cal F}(\Bbb C^N)^{ev}$.
Note that the vaccum vector $\Omega\in {\Cal F}(\Bbb C^N)^{ev}$. 
An operator on ${\Cal F}(\Bbb C^N)$ is called even if it commutes with
$K$. \par
${\Cal F}(\Bbb C^N)$ supports a projective representation of $LSpin(2N)$
at level $1$ (also denoted by $\pi$) , 
and in fact  ${\Cal F}(\Bbb C^N)^{ev}$ is the vacuum representation
of $LSpin(2N)$ (cf. P. 246-7 of [PS]).
\proclaim{Proposition 3.1}
(1):  The vacuum vector $\Omega $ is cyclic and separating for 
$M(I,\Bbb C^N)$ and   $M(I,\Bbb C^N)^\prime = K^{-1} M(I',\Bbb C^N )K$; \par
(2): $M(I,\Bbb C^N)= \pi (L_IU(N))''$; \par
(3): $\pi^{ev} (L_IU(N))'' = \pi(L_ISpin(2N))''$ where
$\pi^{ev} (L_IU(N))''$ denotes the even elements of $\pi (L_IU(N))''$; \par
(4): If   $\Bbb C^N=\Bbb C^{N_1} \oplus \Bbb C^{N_2}$ and let
$\pi(LU(N_1)$ (resp. $\pi(LU(N_2)$) be the representation induced from 
the map
$U(N_1)\rightarrow U(N_1)\oplus id_{N_2}$ (resp. 
$U(N_2)\rightarrow id_{N_1}\oplus U(N_2)$), then
$$
\pi (L_IU(N))''= \pi (L_IU(N_1))'' \vee \pi (L_IU(N_2))''.
$$
\endproclaim
\demo{Proof:}
(1) is proved in \S15 of [W]. (2) is implied in \S15 of [W], 
also cf. Lemma 3.1 of [X6]. To prove (3), note that by (2)
$$
\pi(L_ISpin(2N))'' \subset M^{ev}(I,\Bbb C^N)= \pi^{ev} (L_IU(N))''.
$$ 
Note that both sides are invariant under the action of the modular 
group (cf. [W]), by [T], it is sufficient to show that
$$
\overline{\pi(L_ISpin(2N))''\Omega} \supset  \overline{M^{ev}(I,\Bbb C^N)
\Omega}.
$$
By Reeh-Schlieder theorem 
$$
\overline{\pi(L_ISpin(2N))''\Omega}=\overline{\pi(LSpin(2N))''\Omega}
, $$ and
$$
\overline{\pi(LSpin(2N))''\Omega}=
{\Cal F}(\Bbb C^N)^{ev}
$$ 
by  P. 246-7 of [PS]. Since
$$
\overline{M^{ev}(I,\Bbb C^N)
\Omega} \subset {\Cal F}(\Bbb C^N)^{ev},
$$  the proof of (3) is complete. \par
(4) follows immediately from (2). 
\enddemo \hfill \qed \hfill
\par 
We will consider $N=mn+mk+nk$ and   
$$
\Bbb C^N= \Bbb C^m\otimes \Bbb C^n \oplus \Bbb C^n\otimes \Bbb C^k  
\oplus \Bbb C^m\otimes \Bbb C^k
$$ in the following. \par
Denote by $\pi$ the representation of $LG_1, LG_2$
on ${\Cal F}({\Bbb C}^N)$ induced by the natural inclusions of
$G_1\subset U(N), G_2\subset U(N).$ Note that the levels of representations
match. \par
The $U(1)$ factor of $H_1$ is  mapped into $U(N)$ as
$$
a\rightarrow a^{m+n} id_n \otimes id_m \oplus  a^{n} id_m \otimes id_k 
\oplus a^{-m} id_n \otimes id_k. 
$$
This gives a map $P_1: LU(1)\rightarrow LU(N)$. Denote by
 $P_2: LU(1)\rightarrow LU(N)$ the map induced by
 $$
a\rightarrow a^{m+n} id_m \otimes (id_n \oplus id_k) \oplus id_n\otimes id_k
$$
The  representations $\pi(P_1 (LU(1)))$  and $\pi( P_2(LU(1)))$  have levels
$mn(m+n)(m+n+k)$ and $(m+n)^2 m(n+k)$ respectively. We will denote them by
 $$\pi(LU(1)_{mn(m+n)(m+n+k)})$$ 
and  $\pi(LU(1)_{(m+n)^2 m(n+k)})$ 
 respectively.
We first state a simple result about representations of $LU(1)$:
\proclaim{Lemma 3.2}
If $\pi$ is a positive energy representation of  $LU(1)$, then
it is strongly additive (cf. [L1]), i.e., if  $I_1, I_2$ 
are intervals obtained by removing an interior point
of interval $I$, the 
$$
\pi(L_IU(1))''= \pi(L_{I_1}U(1)''\vee \pi(L_{I_2}U(2))''
$$
\endproclaim
\demo{Proof:}
The representation of the connected component  $LU(1)^0$ of 
 $LU(1)$ is strongly additive by [TL].
Note that  $L_IU(1)$ is generated by  $L_IU(1)^0$ and any loop of winding
number $1$ with support on $I$,  and 
we can choose this loop to have support on $I_1$. This shows that
$$
\pi(L_IU(1))''\subset \pi(L_{I_1}U(1)''\vee \pi(L_{I_2}U(2))''
$$
and completes the proof.
\enddemo \hfill \qed \hfill
\par
\proclaim{Lemma 3.3}
$$
\pi(L_I U(1)_{mn(m+n)(m+n+k)})' \cap \pi(L_I G_1)'' =
\pi(L_I U(1)_{(m+n)^2m(n+k)})' \cap \pi(L_I G_1)''.
$$
\endproclaim
\demo{Proof:}
It is sufficient to show that
$$
\pi(L_I U(1)_{mn(m+n)(m+n+k)})'' \vee \pi(L_I G_1)' =
\pi(L_I U(1)_{(m+n)^2m(n+k)})'' \vee \pi(L_I G_1)'.
$$
Note that for any $\beta\in L_IU(1)$, 
$P_1(\beta)= P_2(\beta) P_3(\beta)$, where $ P_3(\beta): 
P_2(\beta)^{-1}P_1(\beta)\in LU(k)_{m+n}$.
Also $\pi(LU(k)_{m+n})'' \subset  \pi(L_I G_1)'$. Hence
$\pi(P_1(\beta)) \in \pi(P_1(\beta)) \vee  \pi(L_I G_1)'$. This shows
$\subset$ in the lemma. The other inclusion is similar.
\enddemo \hfill \qed \par
\proclaim{Lemma 3.4}
(1). 
$$
\pi(L_ISU(m)_{k+n})' \cap K \pi(L_IU(mn+mk))''K^{-1} =  
K \pi(L_IU(k+n)_m)''K^{-1};
$$ \par
(2).
$$
\align
\pi(L_ISU(m)_{k+n})' & \cap \pi(L_I U(1)_{(m+n)^2m(n+k)})'\cap
K\pi(L_IU(mn+mk))''K^{-1} \\
&=  \pi(L_ISU(k+n))_{m})'';
\endalign
$$ \par
\endproclaim
\demo{Proof:}
Ad (1): Since elements of $\pi(L_ISU(m)_{k+n})''$ commute with $K$, 
it is sufficient to show that:
$$
\pi(L_ISU(m)_{k+n})' \cap  \pi(L_IU(mn+mk))'' =  
\pi(L_IU(k+n)_m)''.
$$ 
By local equivalence (cf. Th. B of [W]), 
it is sufficient to show the above equality for
the restriction $\pi_1$ of $\pi$ to ${\Cal F}({\Bbb C}^{mn+mk})$.
Note that 
$$
\pi_1(L_ISU(m)_{k+n})' \cap  \pi_1(L_IU(mn+mk))'' \supset  
\pi_1(L_IU(n+k)_m)''
$$
and both sides are invariant under the action of modular group. By [T],
it is sufficient to show that
$$
\overline{\pi_1(L_ISU(m)_{k+n})' \cap  \pi_1(L_IU(mn+mk))''\Omega} \subset  
\overline{\pi(L_IU(n+k)_m)''\Omega}
$$
By the decomposition of  ${\Cal F}({\Bbb C}^{mn+mk})$ with respect to
$LSU(m)_{k+n} \times LU(n+k)_m$ given in Prop. 10.6.4 of [PS], 
$\Omega= \Omega_1 \otimes \Omega_2 \in {\Cal H}_1 \otimes {\Cal H}_2,$
where ${\Cal H}_1$ and   ${\Cal H}_2$ are vacuum representations of 
$LSU(m)_{k+n}$ and $LU(n+k)_m$, and $\Omega_1$, $\Omega_2$ are vacuum vectors.
By  Reeh-Schlieder theorem, $
\overline{\pi(L_IU(n+k)_m)''\Omega}= \Omega_1 \otimes {\Cal H}_2.$
Now let $x\in \pi_1(L_ISU(m)_{k+n})' \cap  \pi_1(L_IU(mn+mk))''$, then
$x\in \pi_1(L_ISU(m)_{k+n})' \vee \pi_1(L_{I'}SU(m)_{k+n})' =
\pi_1(LSU(m)_{k+n})'$ by strong additivity (cf. [TL]), and so
$x\Omega \in  \Omega_1 \otimes {\Cal H}_2$, and the proof is complete. \par
Ad (2):
Note that  the right hand side is contained in the 
left  hand side. 
By (1) it is sufficient to show that
$$
\pi(L_I U(1)_{(m+n)^2m(n+k)})' \cap  K \pi(L_IU(k+n)_m)''K^{-1}
\subset \pi(L_ISU(k+n))_{m})''.
$$ 
Note that  both sides are invariant under the action of modular
group. Let $a$ be an element of the  left hand side.
Then  $a \in \pi(LSU(m)_{k+n})' \cap \pi(L U(1)_{(m+n)^2m(n+k)})'$ by 
strong additivity of $LSU(m)$ (cf. [TL]) and 
$LU(1)$ (cf. Lemma 3.2). Now the proof is similar to that of (1).
By using  Reeh-Schlieder theorem and decompositions given in Prop. 10.6.2
and 10.6.4 of [PS], we have  that 
$a\Omega \subset   \overline{\pi(L_ISU(k+n))_{m})''\Omega}.$
By [T], this shows (2).
\enddemo \hfill \qed \par
\proclaim{Lemma 3.5}
$$
\pi(L_IH_1)' \cap \pi(L_IG_1)'' \subset \pi(L_IG_2)''.
$$
\endproclaim
\demo{Proof:}
By Lemma 3.3 
$$
\align
\pi(L_IH_1)' \cap \pi(L_IG_1)'' & \subset \\ 
& \pi(L_ISU(m)_{n+k}' \cap\pi(L_I U(1)_{(m+n)^2m(n+k)})' \cap \pi(L_I G_1)''.
\endalign
$$
Note that by (4) of Prop. 3.1
$$
\align
\pi(L_I G_1)'' & \subset  K \pi(L_IU(mn+mk+nk))'' K^{-1} \\ 
&=  K \pi(L_IU(mn+mk))'' K^{-1} \vee  K \pi(L_IU(nk))'' K^{-1}. 
\endalign
$$
By Lemma 3.4
$$
\align
& \pi(L_IH_1)' \cap \pi(L_IG_1)''  \subset 
\pi(L_ISU(m)_{n+k}' \cap \pi(L_I U(1)_{(m+n)^2m(n+k)})' \cap
\\
& K \pi(L_IU(mn+mk))'' K^{-1} 
\vee  K \pi(L_IU(nk))'' K^{-1} \\
& = \pi(L_ISU(n+k)_{m}'' \vee  K \pi(L_IU(nk))'' K^{-1}.
\endalign
$$
Note that by (3) of Prop. 3.1 
$ \pi(L_ISpin(2nk))''$ are the even elements of $$\pi(L_IU(nk))'',$$ 
and the elements of $\pi(L_IH_1)' \cap \pi(L_IG_1)''$ and 
$\pi(L_ISU(n+k)_{m}''$
are even, 
it follows that
$$
\pi(L_IH_1)' \cap \pi(L_IG_1)'' \subset  \pi(L_ISU(n+k)_{m})'' \vee 
\pi(L_ISpin(2nk))'' = \pi(L_IG_2)''.
$$
\enddemo \hfill \qed \par
\proclaim{Lemma 3.6}
$$
\pi(L_IH_1)' \cap \pi(L_IG_1)'' \subset \pi(L_IH_2)'.
$$
\endproclaim
\demo{Proof:}
By definitions it is enough to show that
$$
\pi(L_IH_1)' \cap \pi(L_IG_1)'' \subset \pi(L_IU(1)_{kn(k+n)(k+n+m)})'
$$ or equivalently
$$
\pi(L_IH_1)'' \vee \pi(L_IG_1)' \supset \pi(L_IU(1)_{kn(k+n)(k+n+m)})''.
$$
Note that $ \pi(L_IU(1)_{kn(k+n)(k+n+m)})$ is actually
$\pi(P(\alpha))$, where $P:LU(1)\rightarrow LU(N)$ is given by
$$
\alpha\rightarrow \alpha^{k+n} id_n \otimes id_k \oplus \alpha^n id_k \otimes
id_m \oplus \alpha^{-k}id_n\otimes id_m.
$$
So 
$$
\align
P(\alpha)&= [\alpha^{k+n} (id_n \oplus id_m) \otimes id_k \oplus 
id_n\times id_m ] \\ 
&\times  
[\alpha^{-k} id_m \otimes (id_k \oplus id_n)\oplus 
id_k\times id_n] 
\endalign
$$
Denote by 
$$
\align
\alpha_1:&= \alpha^{k+n} (id_n \oplus id_m) \otimes id_k \oplus 
id_n\otimes id_m , \\ 
\alpha_2:&=  \alpha^{-k} id_m \otimes (id_k \oplus id_n)\oplus 
id_k\times id_n .
\endalign
$$
Then $\pi(P(\alpha))$ is equal to $\pi(\alpha_1)\pi(\alpha_2) $
up to a scalar. Note that  $$\pi(\alpha_1) \in \pi(LU(k)_{n+m})''$$
and  $\pi(\alpha_2) \in \pi(LU(m)_{k+n})''$, so
$$
\pi(L_IU(1)_{kn(k+n)(k+n+m)})'' \subset 
\pi(LU(k)_{n+m})'' \vee \pi(LU(m)_{k+n})''.
$$
By definition,  $\pi(L_IG_1)'\supset \pi(LU(k)_{n+m})''$, and 
by Lemma 3.5, 
$$
\pi(L_IH_1)'' \vee \pi(L_IG_1)' \supset \pi(L_IG_2)' \supset
\pi(LU(m)_{k+n})''
.
$$
It follows that 
$$
\pi(L_IH_1)'' \vee \pi(L_IG_1)' \supset \pi(L_IU(1)_{kn(k+n)(k+n+m)})''.
$$
\enddemo \hfill \qed \par
\proclaim{Theorem 3.7}
The conformal precosheaves ${\Cal A}(G(m,n,k))$ and  ${\Cal A}(G(k,n,m))$
are isomorphic.
\endproclaim
\demo{Proof:}
By Lemmas 3.5-3.6 ,  for every interval $I$, 
$$
\pi(L_IH_1)' \cap \pi(L_IG_1)'' \subset \pi(L_IH_2)' \cap \pi(L_IG_2)''.
$$
Exchanging $m$ and $k$ in Lemmas 3.5-3.6 , we get
$$
\pi(L_IH_1)' \cap \pi(L_IG_1)'' \supset \pi(L_IH_2)' \cap \pi(L_IG_2)'',
$$
and so
$$
\pi(L_IH_1)' \cap \pi(L_IG_1)'' = \pi(L_IH_2)' \cap \pi(L_IG_2)''.
$$
Let ${\Cal H}$ be the closure of $
\pi(L_IH_1)' \cap \pi(L_IG_1)''\Omega$, and let $P_0$ be the projection
onto  ${\Cal H}$. Let ${\Cal A}$ be the conformal precosheaf given by
$$
{\Cal A}:= \pi(L_IH_1)' \cap \pi(L_IG_1)''P_0 = \pi(L_IH_2)' 
\cap \pi(L_IG_2)''P_0
$$ 
on ${\Cal H}$. 
It follows by definitions that ${\Cal A}(G(m,n,k))$ and  ${\Cal A}(G(k,n,m))$
are both isomorphic to  ${\Cal A}$. 
\enddemo \hfill \qed \par
Note that by Th. 2.4 
${\Cal A}(G(m,n,k))$ has only finitely number of irreducible
representations , and they generate a unitary modular category.
Denote this  modular category by $MC(G(m,n,k))$.
Th. 3.7 implies that  
\proclaim{Corollary 3.8}
There exists a one to one correspondence between the irreducible
representations of  ${\Cal A}(G(m,n,k))$ and  ${\Cal A}(G(k,n,m))$
such that the three manifold invariants (including colored ones, cf. [Tu])
calculated from $MC(G(m,n,k))$ are identical to that from
$MC(G(k,n,m))$. 
\endproclaim
In particular the corollary shows the existence of identifications between all
chiral quantities of  ${\Cal A}(G(m,n,k))$ and  ${\Cal A}(G(k,n,m))$.
By using Th. 4.7 and [X5], one can write down a formula for
the closed three manifold invariants from $MC(G(m,n,k))$. We will
omit the formula, but we note that the symmetry
under the exchange of $m$ and $k$ agrees with \S3 of [X5]. \par
\heading{\S4. Representations of ${\Cal A}(G(m,n,k))$} \endheading
By (2) of Th. 2.4, every irreducible representation of 
${\Cal A}(G(m,n,k))$ occurs in 
$(i,\alpha)$ for some $(i,\alpha)\in exp.$  So we need 
to determine $exp$. It is also known that
there may be field identifications, 
i.e., there may be $(j,\beta)$ with $i\neq j$ or 
$\beta\neq \alpha$ but $(j,\beta)$ is equivalent to 
$(i,\alpha)$ as representations. There are also  issues of fixed point 
resolutions, i.e., as a  representation $(i,\alpha)$ may not
be irreducible, and we need to decompose $(i,\alpha)$ into irreducible
pieces. 
To answer these questions, 
it turns out one needs to determine all VPs for
${\Cal A}(G(m,n,k))$. 
Let us first introduce some notations. \par
Note that the inclusion $H_1\subset G_1$ is
a composition of two inclusions as described at the end of \S1:
 $$   
\align        
H_1  & \subset 
SU(m)_{n} \times SU(m)_{k}  \times SU(n)_{m} \times SU(n)_{k} \\ 
& \times U(1)_{mn(m+n)(m+n)} \times U(1)_{mn(m+n)(k)}  
\endalign
$$ 
and
$$  
\align
&(SU(m)_{n}  \times SU(n)_{m}  
\times U(1)_{mn(m+n)(m+n)})   \times
(SU(m)_{k} \times SU(n)_{k} \\
&  \times U(1)_{mn(m+n)(k)})   
  \subset  G_1:=Spin(2mn)_1 \times SU(m+n)_k. 
\endalign
$$ 
We note that the inclusion
$$
(SU(m)_{n}  \times SU(n)_{m}  
\times U(1)_{mn(m+n)(m+n)}) \subset Spin(2mn)_1
$$
is a conformal inclusion, which is in fact a composition of two
conformal inclusions 
$$
SU(m)_{n}  \times SU(n)_{m}\subset SU(mn)_1
$$
and 
$$
SU(mn)_1\times U(1)_{mn(m+n)(m+n)} \subset Spin(2mn)_1.
$$ 
We will use $\pi_0, \lambda_0, \dot\lambda_1, \dot\lambda_2$ and $\dot q$ to 
denote the representations
of $$Spin(2mn)_1, SU(m+n)_k, SU(m)_{n+k}, SU(n)_{m+k}$$ and 
$ U(1)_{mn(m+n)(m+n+k)}$ respectively. So the general 
coset labels $(i,\alpha)$ in \S2.2 can be identified in the case of 
$H_1\subset G_1$ as $i=\{\pi_0, \lambda_0 \}$ and 
$\alpha= \{\dot\lambda_1, \dot\lambda_2,\dot q \}$. \par
We will use $\lambda_1, \lambda_2, \tilde \lambda_1, \tilde \lambda_2,
\tilde q $ and $q$ to denote the representations of
$$SU(m)_{k}, SU(n)_{k},  SU(m)_{n}, SU(n)_{m}, U(1)_{mn(m+n)(m+n)}$$ 
and  $U(1)_{mnk(m+n)}$ respectively. \par
We use $\tau$ to denote the generator of some symmetries of the 
extended Dynkin diagram of the Kac-Moody algebra, and it is 
defined as follows:\par
Acting on an $SU(K)_M$ representation $\lambda$, $\tau$ rotates the
extended Dynkin indices, i.e., $a_i(\tau(\lambda))= a_{i+1}(\lambda)$,
where $a_{i+K}= a_i$. Acting on the  representations of 
 $Spin(2L)_1$, $\tau$ exchanges the vacuum and vector representations,
and  exchanges the two spinor representations. \par
In accordance with the conventions of \S2.1, if $1$ is used to 
denote a representation of $SU(K)_M$ or $Spin(2L)_1$, it will always
be the vacuum representation. We will however use $0$ to label the
vacuum representation of $U(1)_{2M}$. \par
Denote by $H_3:= SU(m)_{k} \times SU(n)_{k} \times U(1)_{mn(m+n)(k)}$
and $  G_3:=SU(m+n)_k$. \par
\subheading{\S4.1 Selection Rules}
Let $(\lambda_0; \lambda_1,\lambda_2, q)$ be in the $exp$ of $H_3\subset G_3$.
By looking at the actions of the centers of $H_3,G_3$, we can get 
constraints on the
conjugacy classes of the representations. These are known as {\it selection
rules}. For a representation $\lambda$ of  $SU(K)_M$ , we denote by
$r_\lambda$ the number of boxes in the Young tableau corresponding
to $\lambda$. 

\par 
First we have
$$
[e^{\frac{2\pi i}{m}} id_m \oplus id_n]\times
[e^{\frac{-2\pi i n}{m(m+n)}} id_m \oplus e^{\frac{2\pi i}{(m+n)}}id_n] 
= e^{\frac{2\pi i }{(m+n)}} (id_m \oplus id_n).
$$
Note that $e^{\frac{2\pi i}{m}} id_m \oplus id_n$ and 
$e^{\frac{2\pi i }{(m+n)}} (id_m \oplus id_n)$ are in the centers of
$SU(m)$ and $SU(m+n)$ respectively. By considering the actions of these
elements on the space labeled by  $(\lambda_0; \lambda_1,\lambda_2, q),$
we get:
$$
e^{{\frac{2\pi i}{m+n}}r_{\lambda_0}}= e^{\frac{2\pi i}{m}r_{\lambda_1}}
e^{\frac{-2\pi i}{m(m+n)}q},
$$
and so
$$
q =-mr_{\lambda_0} + (m+n)r_{\lambda_1} 
\bmod  m(m+n) \tag 4.1
$$
Similarly by considering the center of $SU(n)$ we get
$$
q  =nr_{\lambda_0}  -(m+n)r_{\lambda_2} 
\bmod  n(m+n) \tag 4.2  
$$
Now let $(\lambda_0, \pi_0; \dot\lambda_1,\dot\lambda_2, \dot q)$ 
be in the set $exp$ of $H_1\subset G_1$.
Similarly by looking at the action of the
centers of $G_1$ and $H_1$ as above , we get the following constraint on the
conjugacy classes of the representations:
$$
\align
\dot q& =-mr_{\lambda_0} + (m+n)r_{\dot\lambda_1} + \frac{1}{2}nm(m+n)\epsilon
\bmod  m(m+n) \tag 4.3 \\
\dot q& =nr_{\lambda_0} - (m+n)r_{\dot\lambda_2} + \frac{1}{2}nm(m+n)\epsilon
\bmod  n(m+n) \tag 4.4  
\endalign
$$
where $\epsilon=1$ or $1$ if $\pi_0$ is a spin representation or otherwise.
\subheading{\S4.2 Vacuum Pairs for $G(m,n,k)$}
Now we are ready to determine the VPs for $H_1\subset G_1$.
By Lemma 2.5, $(\lambda_0,\pi_0;\dot\lambda_1, \dot\lambda_2, \dot q) $ 
is a vacuum pair
iff there exist $\lambda_1, \tilde\lambda_1, \lambda_2, \tilde\lambda_2, 
q,\tilde q$ such that
$(\pi_0;  \tilde\lambda_1,\tilde\lambda_2),$ 
$(\lambda_0; \lambda_1, \lambda_2),$
$(\lambda_1,\tilde\lambda_1;\dot\lambda_1),$
$(\lambda_2,\tilde\lambda_2;\dot\lambda_2)$ 
and $(\tilde q,q;\dot q) $ are VPs. 
By 2.7.12 of [KW], $(\lambda_1,\lambda_2;\dot\lambda_1),$
$(\tilde\lambda_1,\tilde\lambda_2;\dot\lambda_2)$ are VPs iff
$$
(\lambda_1,\tilde\lambda_1;\dot\lambda_1)=(\tau^j(1),\tau^j(1);\tau^j(1))
,(\lambda_2,\tilde\lambda_2;\dot\lambda_2)=
(\tau^i(1),\tau^i(1);\tau^i(1))  
$$
for some $0\leq j\leq m-1, 0\leq i\leq n-1.$
So we should determine VPs of the form $(\tilde q,q;\dot q) $ and 
$(\lambda_0; \tau^j(1), \tau^i(1))$. The following lemma solves the
first question:
\proclaim{Lemma 4.1}
The VPs for the diagonal inclusion 
$$
U(1)_{2a+2b} \subset U(1)_{2a} \times  U(1)_{2b}
$$
are given by
$(\frac{a}{(a,b)} i, \frac{b}{(a,b)} i; \frac{a+b}{(a,b)} i)$
where $ 0\leq i \leq 2(a,b)-1$, and $(a,b)$ is the greatest common
divisor of $a$ and $b$.
\endproclaim
\demo{Proof:}
We use $0\leq x\leq 2a-1, 0\leq y\leq 2b-1, 0\leq z\leq 2(a+b)-1
$ to label the representations. Using 
$$
h_x= \frac{x^2}{4a}, h_y= \frac{y^2}{4b},h_x= \frac{z^2}{4(a+b)} 
$$
one checks easily from (2.4) that the list in the lemma are indeed VPs. We 
want to show that the list is complete. This is an easy exercise and we 
will prove it by calculating  (2.2) in \S2. 
Note that $(x,y;z)\in exp$
if and only if $x+y-z$ is divisible by $2(a,b)$. Note that all the
sectors in this coset have statistical dimensions equal to 1. By lemma
2.2 and (2) of Prop. 3.1 (set $i=1, \alpha=1, z=1$) in [X3], we get
$$
\frac {1}{b(1,1)^2} \times 2(a+b)= 2a \times 2b \times (\frac{a+b}{(a,b)})^2.
$$
It follows that
$$
b(1,1)= \frac{(a,b)}{\sqrt{2ab(a+b)}}.
$$
Note by definition
$$
b(1,1)= \sum_{(x,y;z) \in VPS}\frac{1}{\sqrt{8ab(a+b)}},
$$
and by comparing with the value of $b(1,1)$ we conclude that the
number of VPs must not exceed $2(a,b)$. Thus the list of VPs in the lemma is
complete.
\enddemo \hfill \qed \hfill
\par
In the next few lemmas we determine VPs for $H_3\subset G_3$ of the
form 
$$(\lambda_0; \tau^j(1), \tau^i(1), q).$$ 
Note that the sectors
$ \tau^j(1), \tau^i(1), q$ have statistical dimensions equal to $1$,
and by the argument of Lemma 2.7, such VPs form an abelian group with
group law being the composition of sectors. We denote this abelian
group by ${\Cal S}$.  Also the statistical
dimension of $\lambda_0$ is equal to $1$, so $\lambda_0 x$ must
be irreducible for any sector $x$ of $SU(m+n)_k$. Choose $x$ corresponding
to the fundamental representation of $SU(m+n)$ and using the well known
fusion rules (cf. [W]), we conclude that $\lambda_0= \tau^l(1)$ for
some $ 0\leq l \leq m+n-1.$
We will choose the roots $\alpha_1,..., \alpha_{m+n-1}$ of $SU(m+n)$
such that $\alpha_1,..., \alpha_{m-1}$ and 
$\alpha_{m+1},..., \alpha_{m+n-1}$ 
are roots of $SU(m)$ and $SU(n)$ respectively. We will denote the 
fundamental weights of  $LSU(m+n),  LSU(m)$ and  $LSU(n)$
by $\Lambda_j, \dot\Lambda_{j'}$ and $\ddot\Lambda_{j''}$ respectively, 
where $0\leq j\leq m+n-1, 0\leq j'\leq m-1, 0\leq j''\leq n-1.$
\proclaim{Lemma 4.2}
$(1;\tau^j(1), \tau^i(1), 0) \in {\Cal S}$ iff $j=0 \bmod  m, i=0 \bmod  n.$ 
\endproclaim
\demo{Proof}
Assume that $j\neq 0, i\neq 0$ and that
the coset vacuum vector in ${\bold H}_{(1;\tau^j(1), \tau^i(1), 0)}$
appear in the weight space of $LH_3$ 
with weight $k \Lambda_0 - {\bold r}$, where  
$$
{\bold r} = 
\sum_{0\leq s\leq m+n-1}y_s \alpha_s, y_s \geq 0,0\leq s\leq m+n-1 $$  
Note that $\alpha_0 = \delta - \sum_{1\leq s\leq m+n-1} \alpha_s$ (cf. \S1 of 
[KW]). 
By the  equation (2.4), we get:
$$
\align
(y_0-y_1) \alpha_1 + ...+ (y_0-y_{m-1}) \alpha_{m-1} +  (y_m-y_0) 
\dot\Lambda_{m-1} & = k \dot\Lambda_j \\
(y_0-y_{m+1}) \alpha_{m+1} + ...+ (y_0-y_{m+n-1}) \alpha_{m+n-1} +  (y_m-y_0) 
\ddot\Lambda_{m+1} & = k \ddot\Lambda_i \\ 
y_0-y_m= 0, y_0& = h_{k \dot\Lambda_j} +  h_{k \ddot\Lambda_i}
\endalign
$$
Solving these equations , we get in particular
$$
y_0- y_j = \frac{kj(m-j)}{m}, y_0- y_{m+i} = \frac{ki(n-j)}{n}, 
y_0= \frac{kj(m-j)}{2m} +  \frac{ki(n-j)}{2n}.
$$ 
Note that $ y_j \geq 0,  y_{m+i}\geq 0$ and so
$ y_j=  y_{m+i}=0, y_0=y_m>0.$  It follows that the weight 
 $k \Lambda_0 -{\bold r}
 $ is degenerate (cf. P. 190 of [K]) with respect to 
$k \Lambda_0$, contradicting Lemma 11.2 of [K]. 
\enddemo \hfill \qed \hfill
\par
\proclaim{Lemma 4.3}
If $(\tau^l(1); \tau^j(1),\tau^i(1), q)\in {\Cal S}$, then
$l=j+i \bmod  m+n$ and $q=(nj-mi)k \bmod  mn(m+n)k$.
\endproclaim
\demo{Proof:}
One checks easily using definitions that
$$
(\tau (1); \tau(1), 1, nk)\in {\Cal S}
$$ and 
$$
(\tau (1); 1,\tau(1),  -mk)\in {\Cal S}.
$$ 
Since $S$ is an abelian group, it follows that all
$
(\tau^{j+i} (1); \tau^j(1), \tau^i(1), nj-mi) 
$  form a subgroup ${\Cal S'}$ of ${\Cal S}$. The lemma is equivalent to
${\Cal S'}={\Cal S}.$ 
Without loss of generality let us assume that $n\leq m.$
Let
 $(\tau^l(1); \tau^j(1),\tau^i(1), q)\in {\Cal S},$ to show that
 $(\tau^l(1); \tau^j(1),\tau^i(1), q)\in {\Cal S'},$ by multipying
elements of ${\Cal S'}$ if necessary, we just have to consider the case
$l=0, i=0,$ and we denote by ${\Cal S''}$ the abelian group generated
by such elements.  

Note that $(1; \tau^n(1), 1, n(n+m)k)\in  {\Cal S''}.$ 
Let $(1; \tau^j(1), 1, q)\in  {\Cal S''}$ be an element such that
$q$ is the least positive integer. By Lemma 4.2, ${\Cal S''}$ is a cyclic
group generated by $(1; \tau^j(1), 1, q).$ So there exists a positive
integer $k_1$ such that 
$$  n(n+m)k=qk_1, n=jk_1 \bmod  m.
$$  
To complete the proof we just have to show that $k_1=1$. \par
As in the proof of Lemma 4.2, we have the following  equation  for
$(1; \tau^j(1), 1, q)$ by (2.4):
$$
\align
(y_0-y_1) \alpha_1 + ...+ (y_0-y_{m-1}) \alpha_{m-1} +  (y_m-y_0) 
\dot\Lambda_{m-1}  &= k \dot\Lambda_j \\
(y_0-y_{m+1}) \alpha_{m+1} + ...+ (y_0-y_{m+n-1}) \alpha_{m+n-1} +  (y_m-y_0) 
\ddot\Lambda_{m+1} & = 0 
\endalign
$$
$$ 
y_0-y_m = \frac{q}{m+n}, 
y_0  = h_{k \dot\Lambda_j} +  \frac{q^2}{2mn(m+n)k}.
$$
By solving the equations, we find in particular that
$$
y_0 = \frac{kj(m-j)}{2m} + \frac{q^2}{2(m+n)mnk},
y_0 -y_j = \frac{kj(m-j)}{m} + \frac{qj}{m+n}. 
$$
Since $y_j\geq 0,$ we have the following inequality:
$$
q^2 \geq n(m+n)j(m-j)k^2 + 2nkj q,
$$
and so
$$
q\geq k(nj+\sqrt{njm(m+n-j)}).
$$
Using $ n(n+m)k=qk_1,$ we get
$$
\frac{n(m+n)}{k_1} \geq nj+\sqrt{njm(m+n-j)}.
$$
Solving this equality for $0\leq j \leq m-1,$ we get inequality
$$
j\leq \frac{n}{k_1} + \frac{m}{2} - \sqrt{\frac{m^2}{4}+ \frac{nm(k_1-1)}
{k_1^2}},
$$ 
and so 
$$
k_1j \leq n 
$$ with equality iff $k_1=1.$ Since $k_1j=n \bmod  m$ and $n\leq m$,
we conclude that $k_1=1$.
\enddemo \hfill \qed \par
Now we are ready to prove the following theorem:
\proclaim{Theorem 4.4}
All the VPs of ${\Cal A}(G(m,n,k))$ are given by
$$
(\tau^{j+i}(1), \tau^{jn+im}(1); \tau^{j}(1),\tau^{i}(1),
(nj-mi)(m+n+k))
$$
where $j,i$ are integers.
\endproclaim
\demo{Proof:}
By Lemma 2.5, $(\lambda_0,\pi_0;\dot\lambda_1, \dot\lambda_2, \dot q) $ 
is a vacuum pair
iff there exist 
$$\lambda_1, \tilde\lambda_1, \lambda_2, \tilde\lambda_2, 
q,\tilde q$$
 such that
$(\pi_0;  \tilde\lambda_1,\tilde\lambda_2),$ 
$(\lambda_0; \lambda_1, \lambda_2),$
$(\lambda_1,\tilde\lambda_1;\dot\lambda_1),$
$(\lambda_2,\tilde\lambda_2;\dot\lambda_2)$ 
and $(\tilde q,q;\dot q) $ are VPs. 
By 2.7.12 of [KW], $(\lambda_1,\tilde\lambda_1;\dot\lambda_1),$
$(\lambda_2,\tilde\lambda_2;\dot\lambda_2)$ are VPs iff
$$
(\lambda_1,\tilde\lambda_1;\dot\lambda_1)=(\tau^j(1),\tau^j(1);\tau^j(1))
,(\lambda_2,\tilde\lambda_2;\dot\lambda_2)=
(\tau^i(1),\tau^i(1);\tau^i(1))  
$$
for some $0\leq j\leq m-1, 0\leq i\leq n-1.$
By Lemma 4.1 and Lemma 4.3, we have 
$\lambda_0=\tau^{j+i}(1), \dot\lambda_1=\tau^{j}(1) , \dot\lambda_2
=\tau^{i}(1)$ and $q=(nj-mi)(m+n+k).$ Since
$(\pi_0;  \tilde\lambda_1,\tilde\lambda_2)$ is the VP associated to
a regular conformal inclusion, it is determined by Prop. 4.2 of [KW], and one
checks easily that $\pi_0$ takes the form stated in the theorem.
\enddemo \hfill \qed \par
\proclaim{Corollary 4.5}
(1):
Assume that $(\lambda_0, \pi_0; \dot\lambda_1, \dot\lambda_2, \dot q)$
verifies selection rules (4.3) and (4.4). Then  
$$
b(\lambda_0, \pi_0; \dot\lambda_1, \dot\lambda_2, \dot q)   
=d_{\lambda_0}  d_{\dot\lambda_1} d_{\dot\lambda_2}
b(1,1;1,1,0);
$$\par
(2):  $(\lambda_0, \pi_0; \dot\lambda_1, \dot\lambda_2, \dot q)\in exp$
if and only if it verifies selection rules (4.3) and (4.4). \par
(3) The statistical dimension of $
(\lambda_0, \pi_0; \dot\lambda_1, \dot\lambda_2, \dot q)$
is $d_{\lambda_0}  d_{\dot\lambda_1} d_{\dot\lambda_2}.
$
\endproclaim
\demo{Proof:}
Ad (1): 
To save some writing denote by $i:=\{\lambda_0, \pi_0 \},
\alpha:=  \{ \dot\lambda_1, \dot\lambda_2, \dot q \}.$
By definition   
$$
b(\lambda_0, \pi_0; \dot\lambda_1, \dot\lambda_2, \dot q)= 
\sum_{w\in VPS}S_{i w(1)}\overline{ \dot S_{\alpha w(1)}}. 
$$
Using Th. 4.4,  the assumption and symmetry properties of $S$ matrices
(cf. \S2 of [NS]),
we conclude that
$$
b(\lambda_0, \pi_0; \dot\lambda_1, \dot\lambda_2, \dot q)   
=d_{\lambda_0}  d_{\dot\lambda_1} d_{\dot\lambda_2}
b(1,1;1,1,0).
$$
\par
Ad (2): This follows immediately from (1) and Th. B of [KW].
\par
Ad (3): This follows from (1) and (3) of Th. 2.4 .
\enddemo \hfill \qed 
\par
Assume that $(\lambda_0, \pi_0; \dot\lambda_1, \dot\lambda_2, \dot q)$
verifies selection rules (3) and (4). By Cor. 4.5, 
 $(\lambda_0, \pi_0; \dot\lambda_1, \dot\lambda_2, \dot q)\in exp.$
We will determine the irreducible components of representation
$$(\lambda_0, \pi_0; \dot\lambda_1, \dot\lambda_2, \dot q).$$
To save some writing denote by $i:=\{\lambda_0, \pi_0 \},
\alpha:=  \{ \dot\lambda_1, \dot\lambda_2, \dot q \}.$
By (2) of Prop. 2.3, $a_{(i,\alpha)\otimes 1}\prec 
a_{1\otimes \bar\alpha}\sigma_i,$ but by (3) of Cor. 4.5, 
$d_{(i,\alpha)}= d_i d_\alpha,$ it follows that
$$
a_{(i,\alpha)\otimes 1}=a_{1\otimes \bar\alpha}\sigma_i.
$$ 
So by the same argument as in the derivation of (**) in [X1] and 
use Th. 4.4 we get
$$
\align
\langle (i,\alpha), (i',\alpha') \rangle 
&= \langle a_{1\otimes \bar\alpha}\sigma_i,  
a_{1\otimes \bar\alpha'}\sigma_{i'} \rangle  \\
&= \sum_{w\in VPS} \delta_{w(i),i'} \delta_{w(\alpha),\alpha'} \tag 4.5
\endalign
$$
By setting $i=i', \alpha=\alpha'$ in (4.5), we get
$$
\langle (i,\alpha), (i,\alpha) \rangle 
= t
$$ where $t$ is the number of elements in the set 
$$
F(i,\alpha):= \{ w\in VPS, w(i)=i, w(\alpha)=(\alpha) \} \tag 4.6
$$
\proclaim{Lemma 4.6}
$F(i,\alpha)$ is a cyclic group of order $t$. Moreover, let
$\{j,\beta \}$ be the generator. Then $\sigma_j= a_{1\otimes \beta}$ 
has order $t$,
i.e., $t$ is the least positive integer such that
$\sigma_j^t=1.$ 
\endproclaim
\demo{Proof:}
Let $w\in F(i,\alpha)$.
By Th. 4.4 and  definitions one checks easily the following property:\par
If $\{1, \pi_0; \dot\lambda_1,  \dot\lambda_2, \dot q \} \in F(i,\alpha),$ then $\pi_0=1, \dot\lambda_1=1, \dot\lambda_2=1,  \dot q=0 $; \par
It follows from  that the projection of $w\in F(i,\alpha)$ onto its
first component in ${\Bbb Z}_{m+n}$ is an embedding, and so
$ F(i,\alpha)$ must be a cyclic group of order $t$ which is a divisor
of $m+n$. \par
Now let $\{j,\beta \}$ be the order $t$ generator of $F(i,\alpha).$ 
So $j^t=1, \beta^t=1$ and $t$ is the least positive 
integer with this property.
Let $t_1, t_2$ be the orders of $j, \beta$ respectively. Then $t$ is
the least common multiple of $t_1$ and $t_2$. Since 
$\{j,\beta \} \in VPS,$
$\sigma_j = a_{1\otimes \bar\beta},$ and it follows that
$t_1$ is a divisor of $t_2$ since $j\rightarrow \sigma_j$ is 
an embedding. Note that
$$
a_{1\otimes \bar\beta^{t_1}} = \sigma_j^{t_1} =1,
$$
and so $\{1,\beta^{t_1} \} \in VPS$ and also fix $(i,\alpha)$, 
by the property  above we must have $\beta^{t_1}=1$, so 
$t_2$ is also a  divisor of $t_1.$ It follows that $t_1=t_2=t$.
\enddemo \hfill \qed 
\par
By the formula before  (4.5) the map
$$
(i,\alpha)\rightarrow  a_{1\otimes \bar\alpha} \sigma_i
$$
is a ring isomorphism. 
By definitions, $\sigma_i\sigma_j= \sigma_i$ and 
$$a_{1\otimes \bar\alpha}  \sigma_j =
a_{1\otimes \bar\alpha} a_{1\otimes \bar\beta}
=a_{1\otimes \bar\alpha} 
$$
where $\{j,\beta \}$ is as in Lemma 4.6. 
Moreover, by (4.5) 
$$
\langle a_{1\otimes \bar\alpha}\sigma_i,  a_{1\otimes \bar\alpha}\sigma_i
\rangle =t
$$
and by Lemma 4.6 $\sigma_j$ has order $t$.  
Appling Lemma 2.1 of [X2] in the present case with  
$a= \sigma_i, b=a_{1\otimes \bar\alpha}$
and $\tau=\sigma_j$, we conclude that  
the representation $(i,\alpha)$ decomposes into
$t$ distinct irreducible pieces and each irreducible piece has equal
statistical dimension. We record this result in the following:
\proclaim{Theorem 4.7}
Let  $(\lambda_0, \pi_0; \dot\lambda_1, \dot\lambda_2, \dot q)\in exp $ 
be a representation of ${\Cal A}(G(m,n,k))$. Then the representation 
decomposes into $t$
distinct  
irreducible representations  where $t$ is the number of elements in the 
set (4.6), and each such irreducible  
representation has equal
statistical dimension.
\endproclaim
We note that by (2) of Th. 2.4, Th. 4.7 and formula (4.5) above,
we can  give a list of  all the irreducible
representations of ${\Cal A}(G(m,n,k))$ as follows: \par
First we write down all  
$$(\lambda_0, \pi_0; \dot\lambda_1, \dot\lambda_2, \dot q)$$
which verifies (4.3) and (4.4). Denote such a set by $exp$. $exp$
admits a natural action of $VPS$ given in Th. 4.4. Suppose that $exp$
is the union of $l$ orbits $exp_1,...,exp_l$. Let $(i_p,\alpha_p)\in
exp_p, 1\leq p\leq l$ be representatives of the orbits. 
We note that two different representative of the same orbit are unitarily
equivalent  representations of  ${\Cal A}(G(m,n,k))$ by (4.5)
and Th. 4.7.
Let $t_p$ be
the order of $F(i_p,\alpha_p),  1\leq p\leq l$ as defined in (4.6). 
Then each representation $(i_p,\alpha_p)$ of  ${\Cal A}(G(m,n,k))$
decomposes into  $t_p$ distinct  irreducible pieces, and hence the number
of irreducible  representations of  ${\Cal A}(G(m,n,k))$ is given by
$
\sum_{1\leq p\leq l} t_p.$

\heading References \endheading
\roster

\item"[DJ]" D. Dunbar and K. Joshi,
{\it Characters for coset conformal field theories and Maverick
examples}, Inter. J. Mod. Phys. A, Vol.8, No. 23 (1993), 4103-4121.
\item"{[FZ]}" I. Frenkel and Y. Zhu, 
{\it Vertex operator algebras associated to representations of 
affine and Virasoro algebras}, Duke Math. Journal (1992), Vol. 66, No. 1
, 123-168.
\item"{[FS]}" J. Fuchs and C. Schweigert, {\it 
Level-rank duality of WZW theories and 
isomorphisms of $N=2$ coset models,} Ann.
Physics 234 (1994), no. 1, 102--140.
\item"{[GL1]}"  D.Guido and R.Longo, {\it  The Conformal Spin and
Statistics Theorem},  \par
Comm.Math.Phys., {\bf 181}, 11-35 (1996) 
\item"{[GL2]}"  D.Guido and R. Longo,{\it Relativistic invariance
and charge conjugation in quantum field theory},
Comm.Math.Phys., {\bf 148}, 521-551 (1992).

\item"{[GKO]}" P. Goddard and D. Olive, eds., { \it Kac-Moody and
Virasoro algebras,} Advanced Series in Math. Phys., Vol 3, World
Scientific 1988.
\item"{[H]}" R. Haag, {\it Local Quantum Physics}, Springer-Verlag 1992.
\item"{[K]}"  V. G. Kac, {\it Infinite Dimensional Lie Algebras}, 3rd
Edition,     
Cambridge University Press, 1990.
\item"{[KW]}"  V. G. Kac and M. Wakimoto, {\it Modular and conformal
invariance constraints in representation theory of affine algebras},  
Advances in Math., {\bf 70}, 156-234 (1988).
\item"{[KLM]}" Y. Kawahigashi, R. Longo and M. M\"{u}ger,
{\it Multi-interval Subfactors and Modularity of Representations in
Conformal Field theory}, to appear \par
in Comm.Math.Phys., also see math.OA/9903104.
\item"{[KS]}"  Y. Kazama and H. Suzuki, {\it  
New $N=2$ superconformal field theories and 
superstring compactification,} Nuclear Phys. B
321 (1989), no. 1, 232--268. 
\item"{[LVW]}" W. Lerche, C. Vafa and N. P. Warner, Nucl. Phys. B324
(1989) 427.
\item"{[L1]}"  R. Longo, 
{\it  Conformal Subnets and Intermediate Subfactors,}\par
math.OA/0102196.
\item"{[L2]}"  R. Longo, {\it Duality for Hopf algebras and for
subfactors}, 
I, Comm. Math. Phys., {\bf 159}, 133-150 (1994).
\item"{[L3]}"  R. Longo, {\it Index of subfactors and statistics of
quantum fields}, I, Comm. Math. Phys., {\bf 126}, 217-247 (1989.
\item"{[L4]}"  R. Longo, {\it Index of subfactors and statistics of
quantum fields}, II, Comm. Math. Phys., {\bf 130}, 285-309 (1990).

\item"{[L5]}"  R. Longo, {\it Minimal index and braided subfactors,
} J.Funct.Analysis {\bf 109} (1992), 98-112.
\item"{[LR]}"  R. Longo and K.-H. Rehren, {\it Nets of subfactors},
Rev. Math. Phys., {\bf 7}, 567-597 (1995).
\item"[NS]" S. Naculich and H. Schnitzer, {\it Superconformal 
coset equivalence from level-rank duality}, Nuclear Phys. B 505
(1997), no. 3, 727--748. 
\item"{[PP]}" M.Pimsner and S.Popa,
{\it Entropy and index for subfactors}, \par
Ann. Sci.\'{E}c.Norm.Sup. {\bf 19},
57-106 (1986). 
\item"[PS]" A. Pressley and G. Segal, {\it Loop Groups,} O.U.P. 1986.
\item"[T]"  M. Takesaki, {\it Conditional expectation in von Neumann
algebra,} J. Funct. Analysis 9 (1972), 306-321.

\item"{[TL]}" V. Toledano Laredo, {\it Fusion of
Positive Energy Representations of $LSpin_{2n}$}.
Ph.D. dissertation, University of Cambridge, 1997
\item"[Tu]" V. G. Turaev, {\it Quantum invariants of knots and
3-manifolds,} Walter de Gruyter, Berlin, New York 1994.

\item"{[W]}"  A. Wassermann, {\it Operator algebras and Conformal
field theories III},  Invent. Math. Vol. 133, 467-539 (1998)

\item"[X1]" F.Xu, {\it Algebraic coset conformal field theories},
Comm. Math. Phys. 211 (2000) 1-43.
\item"[X2]" F.Xu, {\it Algebraic coset conformal field theories II},
Publ. RIMS, vol. 35 (1999), 795-824.
\item"[X3]" F.Xu, {\it On a conjecture of Kac-Wakimoto},
Publ. RIMS, vol. 37 (2001), 165-190.
\item"[X4]" F.Xu, {\it   New braided endomorphisms from conformal
inclusions, } \par
Comm.Math.Phys. 192 (1998) 349-403.

\item"[X5]" F.Xu, {\it 3-manifold invariants from cosets},
math.GT/9907077.
\item"[X6]" F.Xu, {\it Jones-Wassermann subfactors for 
Disconnected Intervals}, \par
Comm. Contemp. Math. Vol. 2, No. 3 (2000) 307-347.

\endroster 
\enddocument